\documentclass[11pt]{article}
\usepackage[margin=1in]{geometry}
\usepackage{amsmath,amssymb,amsthm,mathtools}
\usepackage{booktabs,microtype}
\usepackage[colorlinks=true,linkcolor=blue,citecolor=blue,urlcolor=blue]{hyperref}

\newtheorem{theorem}{Theorem}[section]
\newtheorem{proposition}[theorem]{Proposition}
\newtheorem{lemma}[theorem]{Lemma}
\newtheorem{corollary}[theorem]{Corollary}
\theoremstyle{definition}
\newtheorem{definition}[theorem]{Definition}
\newtheorem{example}[theorem]{Example}
\theoremstyle{remark}
\newtheorem{remark}[theorem]{Remark}

\newcommand{\E}{\mathbb E}
\newcommand{\Var}{\operatorname{Var}}

\newcommand{\Comp}{\mathcal C}
\newcommand{\Prob}{\mathbb P}
\newcommand{\N}{\mathcal N}
\newcommand{\ind}{\mathbf 1}
\newcommand{\op}{\mathcal T}

\title{The Gini Index of an Integer Composition:\\
Rational Generating Functions and Limit Theorems}
\author{Grant Kopitzke\\
\small Department of Mathematics, University of Wisconsin--Stevens Point\\
\small Stevens Point, Wisconsin 54481, USA\\
\small \href{mailto:gkopitzk@uwsp.edu}{gkopitzk@uwsp.edu}}
\date{}

\begin{document}
\maketitle

\begin{abstract}
We extend the discrete Gini statistic on integer partitions to integer
compositions by rearranging the parts in decreasing order. Under the uniform
distribution on the $2^{n-1}$ compositions of $n$, the resulting statistic has
a rational ordinary generating function for its aggregate values, as well
as a rational generating function for its second moments. More generally,
the generating function of every fixed raw moment is rational, although
the bivariate generating function encoding the entire distribution is not.
We obtain an exact formula for the mean and prove that the variance of the
unnormalized statistic is asymptotic to $10n^3/189$. The normalization used
for partitions has mean tending to $2/3$, in contrast to the limit $1$ under
uniform random partitions. Using the separator representation of a random
composition and a renewal-reward argument, we prove a central limit theorem:
its centered, normalized Gini statistic has asymptotic Gaussian fluctuations
with variance constant $40/189$. We also prove that the ordinary classical
Gini coefficient of the composition parts tends to $1/3$ and satisfies a
central limit theorem with variance constant $20/189$.
\end{abstract}


\section{Introduction}\label{sec:intro}

An \emph{integer composition} of $n\geq1$ is an ordered sequence
$\alpha=(\alpha_1,\ldots,\alpha_k)$ of positive integers with
$\alpha_1+\cdots+\alpha_k=n$; we write $\alpha\models n$ and denote
its number of parts by $K(\alpha)=k$. Unlike a partition, the order of
its parts matters. For example, $(3,1)$ and $(1,3)$ are different
compositions of $4$, whereas both correspond to the same integer partition.
There are $2^{n-1}$ compositions of $n$, because any subset of the $n-1$
gaps between consecutive units can be chosen as separators
\cite{AndrewsEriksson,DLMFCompositions}.

In \cite{Kopitzke2020}, a discrete Gini statistic was defined for a
partition $\lambda=(\lambda_1,\ldots,\lambda_k)\vdash n$, with
$\lambda_1\geq\cdots\geq\lambda_k$, by
\begin{equation}\label{eq:partitionGini}
 g(\lambda)=\binom{n+1}{2}-\sum_{i=1}^k i\lambda_i.
\end{equation}
It takes values between $0$ and $\binom n2$. Its normalized version,
$g(\lambda)/\binom n2$, uses a normalization depending on the total
size $n$, not the number $k$ of parts. This distinction matters when
comparing it to the conventional Gini coefficient of a sample.

Let $\alpha^\downarrow$ be the integer partition obtained by sorting the
parts of a composition $\alpha$ into weakly decreasing order. The subject
of this paper is the statistic
\begin{equation}\label{eq:compGini}
 g_C(\alpha):=g(\alpha^\downarrow),\qquad \alpha\models n.
\end{equation}
In particular, permuting the parts of a composition does not change
$g_C$. Nevertheless, sampling compositions uniformly is not equivalent
to sampling partitions uniformly: a partition with $k$ parts and
multiplicities $m_1,m_2,\ldots$ has $k!/\prod_jm_j!$ distinct
rearrangements. Uniform composition sampling therefore puts more weight
on partitions having many rearrangements.

This change of probability space has substantial consequences. Under
uniform sampling of partitions of $n$, the normalized statistic from
\eqref{eq:partitionGini} tends to $1$ in probability, whereas the
corresponding composition statistic tends to $2/3$. The composition
model also admits rational generating functions for all fixed raw
moments. This allows exact and asymptotic moment calculations by
ordinary coefficient extraction, without invoking the saddle-point
methods required for Euler's partition product.

The main results are the first- and second-moment generating functions,
Theorems~\ref{thm:firstGF} and~\ref{thm:secondGF}, the asymptotic moment
formulas in Theorem~\ref{thm:moments}, and the Gaussian limit theorem
in Theorem~\ref{thm:jointCLT}. The last theorem also yields a central
limit theorem for the conventional Gini coefficient of the composition
parts. Random compositions and their geometric structure have been
investigated in various settings; see, for instance,
\cite{BevanThrelfall2025}. Our use of geometric run lengths makes it
possible to prove the Gaussian fluctuation results directly.

\section{Two descriptions of the composition Gini statistic}\label{sec:identities}

We begin with an identity that removes the apparent difficulty caused
by sorting the parts of a composition. It is also useful for separating
this statistic from the conventional sample Gini coefficient.

\begin{proposition}\label{prop:pairwise}
Let $\alpha=(\alpha_1,\ldots,\alpha_k)\models n$, and define
\[
 D(\alpha)=\sum_{1\leq i<j\leq k}|\alpha_i-\alpha_j|.
\]
Then
\begin{equation}\label{eq:pairwise}
 g_C(\alpha)=\frac{n(n-k)+D(\alpha)}2.
\end{equation}
\end{proposition}
\begin{proof}
Write $a_i=\alpha_i^\downarrow$. Since the $a_i$ are decreasing,
$|a_i-a_j|=a_i-a_j$ whenever $i<j$. In the sum over $i<j$,
the part $a_i$ occurs positively $k-i$ times and negatively $i-1$
times. Consequently,
\[
 D(\alpha)=\sum_{i=1}^k(k+1-2i)a_i
 =(k+1)n-2\sum_{i=1}^k ia_i.
\]
Solving for the final sum and inserting the result into
\eqref{eq:partitionGini} gives
\[
 g_C(\alpha)
 =\binom{n+1}{2}-\frac{(k+1)n-D(\alpha)}2
 =\frac{n(n-k)+D(\alpha)}2.
\]
\end{proof}

\begin{definition}\label{def:classicalGini}
For a composition $\alpha\models n$ with $k$ parts, define its
\emph{classical (uncorrected) Gini coefficient} by
\begin{equation}\label{eq:classicalGini}
 G_{\mathrm{cl}}(\alpha)=\frac{D(\alpha)}{kn}.
\end{equation}
For $n\geq2$ we also define the normalized composition statistic
\begin{equation}\label{eq:normalizedGini}
 \widetilde g_C(\alpha)=\frac{g_C(\alpha)}{\binom n2}.
\end{equation}
\end{definition}

The normalization in \eqref{eq:classicalGini} is the usual
pairwise-absolute-difference normalization for $k$ observations with
total $n$. Combining the definitions with Proposition~\ref{prop:pairwise}
yields
\begin{equation}\label{eq:twoGinis}
 \widetilde g_C(\alpha)
 =\frac{n-k}{n-1}+\frac{k}{n-1}G_{\mathrm{cl}}(\alpha).
\end{equation}
Thus the two notions of Gini coefficient are related, but generally
not equal. In particular, if all $k$ parts are equal, then
$G_{\mathrm{cl}}=0$, while $\widetilde g_C=(n-k)/(n-1)$.

\begin{example}\label{ex:small}
The eight compositions of $4$ have the following values:
\[
\begin{array}{c|c|c}
 \alpha & \alpha^\downarrow & g_C(\alpha)\\ \hline
 (4)&(4)&6\\
 (3,1),(1,3)&(3,1)&5\\
 (2,2)&(2,2)&4\\
 (2,1,1),(1,2,1),(1,1,2)&(2,1,1)&3\\
 (1,1,1,1)&(1,1,1,1)&0
\end{array}
\]
The sum of the eight values is $29$, and the sum of their squares
is $129$. Thus the mean is $29/8$ and the variance is $191/64$,
providing small checks on the first two generating functions below.
The multiplicities of the rearrangements explain why this is not the
same as averaging the five partition values with equal weights.
\end{example}

There is a second description in terms of the columns of a Young
diagram. For $j\geq1$, let
\[
 c_j(\alpha)=\#\{i:\alpha_i\geq j\}.
\]
These are the parts of the conjugate of $\alpha^\downarrow$.
The row--column identity
$\sum_i ia_i=\sum_j\binom{c_j+1}{2}$ and
$\sum_jc_j=n$ give
\begin{equation}\label{eq:columnGini}
 g_C(\alpha)
 =\binom{n+1}{2}-\sum_{j\geq1}\binom{c_j+1}{2}
 =\frac12\left(n^2-\sum_{j\geq1}c_j^2\right).
\end{equation}
Equivalently, $g_C(\alpha)=e_2(c_1,c_2,\ldots)$.
We shall return to \eqref{eq:columnGini} in Section~\ref{sec:CLT}:
it converts the asymptotic problem into the fluctuations of a
collection of part-counting statistics.

\section{The distribution and the first-moment generating function}
\label{sec:first}

Let $\Comp_n=\{\alpha:\alpha\models n\}$. The full bivariate
generating function is
\begin{equation}\label{eq:fullDistribution}
 \mathcal H(q,x)=1+\sum_{n\geq1}\sum_{\alpha\in\Comp_n}
      q^{g_C(\alpha)}x^n.
\end{equation}
One may regroup this series by the associated sorted partition:
\begin{equation}\label{eq:regroup}
 \mathcal H(q,x)=1+\sum_{n\geq1}x^n
 \sum_{\lambda\vdash n}\frac{\ell(\lambda)!}
 {\prod_{j\geq1}m_j(\lambda)!}\,q^{g(\lambda)},
\end{equation}
where $m_j(\lambda)$ counts parts equal to $j$. For instance,
\[
 [x^4]\mathcal H(q,x)=1+3q^3+q^4+2q^5+q^6.
\]
The bivariate series retains the entire distribution, but the
ordinary series obtained by weighting the compositions by $g_C$
turns out to be much simpler.

\begin{theorem}\label{thm:firstGF}
Let
\[
 F_1(x)=\sum_{n\geq1}\left(\sum_{\alpha\models n}
                         g_C(\alpha)\right)x^n.
\]
Then
\begin{equation}\label{eq:F1}
 \boxed{F_1(x)=\frac{x^2(1+2x)}{(1+x)(1-2x)^3}.}
\end{equation}
\end{theorem}

\begin{proof}
A positive part contributes
\[
 A(x)=x+x^2+\cdots=\frac{x}{1-x}.
\]
As a composition is an ordered sequence of positive parts, the
generating function for compositions, including the empty one, is
\begin{equation}\label{eq:Cx}
 C(x)=1+A(x)+A(x)^2+\cdots
     =\frac1{1-A(x)}=\frac{1-x}{1-2x}.
\end{equation}
More generally, marking the number $k$ of parts by $u$ gives
\begin{equation}\label{eq:Cxu}
 C(x,u)=\frac1{1-uA(x)}
        =\frac{1-x}{1-(1+u)x}.
\end{equation}
Introduce $\Theta_x=x\partial/\partial x$ and
$\Theta_u=u\partial/\partial u$. Since
\[
 (\Theta_x^2-\Theta_x\Theta_u)(x^nu^k)
       =n(n-k)x^nu^k,
\]
the generating function for the contribution $n(n-k)/2$
in \eqref{eq:pairwise} is
\begin{equation}\label{eq:partJ}
 \left.\frac12(\Theta_x^2-\Theta_x\Theta_u)C(x,u)
 \right|_{u=1}=\frac{x^2}{(1-2x)^3}.
\end{equation}

We next enumerate the absolute-difference contribution $D(\alpha)$.
For two distinguished positive parts, set
\[
 B(x)=\sum_{a,b\geq1}|a-b|x^{a+b}.
\]
The terms with $a=b$ vanish. By symmetry, put $b=a+h$, $h\geq1$:
\begin{align*}
 B(x)
 &=2\sum_{a,h\geq1}h x^{2a+h}
 =2\left(\sum_{a\geq1}x^{2a}\right)
     \left(\sum_{h\geq1}hx^h\right)\\
 &=\frac{2x^3}{(1-x)^2(1-x^2)}.
\end{align*}
Now choose two positions $i<j$ in a composition. There is an
arbitrary (possibly empty) composition before the first position,
another between the two positions, and a third after the second.
These three unrestricted portions contribute $C(x)^3$, while the
marked positions contribute $B(x)$. Summing over all choices of
the positions therefore gives
\begin{equation}\label{eq:FD}
 \sum_{n\geq1}\sum_{\alpha\models n}D(\alpha)x^n
 =B(x)C(x)^3=\frac{2x^3}{(1+x)(1-2x)^3}.
\end{equation}
Equation~\eqref{eq:pairwise} and
\eqref{eq:partJ}--\eqref{eq:FD} now give
\[
 F_1(x)=\frac{x^2}{(1-2x)^3}
       +\frac{x^3}{(1+x)(1-2x)^3},
\]
which simplifies to \eqref{eq:F1}.
\end{proof}

\begin{corollary}\label{cor:exactMean}
Choose $C_n$ uniformly from $\Comp_n$, and write
$X_n=g_C(C_n)$. For every $n\geq1$,
\begin{equation}\label{eq:exactMean}
 \boxed{\E[X_n]=\frac{n^2}{3}-\frac{4n}{9}
   +\frac{2}{27}-\frac{2(-1)^n}{27\cdot2^n}.}
\end{equation}
\end{corollary}
\begin{proof}
Each of the $2^{n-1}$ compositions has probability $2^{-(n-1)}$,
so $\E[X_n]=[x^n]F_1(x)/2^{n-1}$. The partial-fraction expansion
\[
 F_1(x)=\frac{1}{3(1-2x)^3}
 -\frac{13}{18(1-2x)^2}
 +\frac{23}{54(1-2x)}-\frac1{27(1+x)}
\]
and the coefficient identities
\[
 [x^n](1-2x)^{-r}=2^n\binom{n+r-1}{r-1},
 \qquad [x^n](1+x)^{-1}=(-1)^n
\]
give the asserted formula after simplification.
\end{proof}

\section{Second moments and rationality of all fixed moments}
\label{sec:second}

The first moment requires only one selected pair of positions.
For the second moment, two selected pairs may coincide, share one
position, or have four distinct positions. Separating these cases
leads to another rational generating function.

\begin{theorem}\label{thm:secondGF}
Let
\[
 F_2(x)=\sum_{n\geq1}\left(\sum_{\alpha\models n}
                         g_C(\alpha)^2\right)x^n.
\]
Then
\begin{equation}\label{eq:F2}
 \boxed{
 F_2(x)=\frac{x^2R(x)}{(1+x)^3(1+x+x^2)(1-2x)^5},
 }
\end{equation}
where
\begin{equation}\label{eq:Rpoly}
 R(x)=1+11x+34x^2+71x^3+83x^4+66x^5+20x^6+8x^7.
\end{equation}
\end{theorem}

\begin{proof}
Write $J(\alpha)=n(n-K(\alpha))$ and $D=D(\alpha)$.
Equation~\eqref{eq:pairwise} implies
\begin{equation}\label{eq:squareDecomp}
 4g_C(\alpha)^2=J(\alpha)^2+2J(\alpha)D(\alpha)
                         +D(\alpha)^2.
\end{equation}
For brevity, set $\op=\Theta_x^2-\Theta_x\Theta_u$.
The generating function for $J^2$ is
\begin{equation}\label{eq:J2gf}
 \left.\op^2C(x,u)\right|_{u=1}
 =\frac{2x^2(2x^2+7x+2)}{(1-2x)^5}.
\end{equation}
To treat the mixed term, use $u$ to mark the number of parts
in the two-distinguished-position construction. The corresponding
series is
\[
 \mathcal D(x,u)=u^2B(x)C(x,u)^3.
\]
Applying $\op$ marks $n(n-K)$ and gives
\begin{equation}\label{eq:JDgf}
 \left.\op\mathcal D(x,u)\right|_{u=1}
 =\frac{2x^3(6x^4+7x^3+23x^2+7x+3)}
        {(1+x)^3(1-2x)^5}.
\end{equation}

It remains to calculate the contribution from $D^2$. First define
\begin{align*}
 B_2(x)&=\sum_{a,b\geq1}|a-b|^2x^{a+b},\\
 H_3(x)&=\sum_{a,b,c\geq1}|a-b|\,|a-c|x^{a+b+c}.
\end{align*}
The same difference substitution used for $B(x)$ gives
\begin{equation}\label{eq:B2}
 B_2(x)=2\frac{x^2}{1-x^2}\sum_{h\geq1}h^2x^h
       =\frac{2x^3}{(1-x)^4}.
\end{equation}
For $H_3$, temporarily fix $a$ and introduce
\[
 V_a(x)=\sum_{b\geq1}|a-b|x^b.
\]
Using $|a-b|=(a-b)+2(b-a)\ind_{\{b>a\}}$ and summing
geometric series,
\begin{equation}\label{eq:Va}
 V_a(x)=\frac{ax}{1-x}-\frac{x}{(1-x)^2}
                         +\frac{2x^{a+1}}{(1-x)^2}.
\end{equation}
Consequently, $H_3(x)=\sum_{a\geq1}x^aV_a(x)^2$.
Expanding the square and applying the elementary identities
\[
 \sum_{a\geq1}x^a=\frac{x}{1-x},\quad
 \sum_{a\geq1}ax^a=\frac{x}{(1-x)^2},\quad
 \sum_{a\geq1}a^2x^a=\frac{x(1+x)}{(1-x)^3}
\]
with $x$ replaced by $x^2$ and $x^3$ where needed gives
\begin{equation}\label{eq:H3}
 H_3(x)=\frac{x^4(1+3x+8x^2+3x^3+x^4)}
             {(1-x)^5(1+x)^2(1+x+x^2)}.
\end{equation}

In $D(\alpha)^2$, the two ordered pairs of positions can be
identical; they can be distinct but share one position; or they
can be disjoint. In the first case there is one selected pair and
three intervening compositions, giving $B_2C^3$.
In the second case there are three selected positions and six
ordered choices of two distinct pairs, giving $6H_3C^4$.
In the third case there are four selected positions and six ordered
ways of pairing them, giving $6B^2C^5$. Thus
\begin{equation}\label{eq:D2gf}
 \sum_{n\geq1}\sum_{\alpha\models n}D(\alpha)^2x^n
 =B_2(x)C(x)^3+6H_3(x)C(x)^4+6B(x)^2C(x)^5.
\end{equation}
The powers of $C(x)$ count the unrestricted compositions between,
before, and after the selected positions.
Finally, substitute \eqref{eq:J2gf}, \eqref{eq:JDgf}, and
\eqref{eq:D2gf} into \eqref{eq:squareDecomp}, divide by $4$,
and simplify the resulting rational expression. The numerator is
$x^2R(x)$, and the denominator is
$(1+x)^3(1+x+x^2)(1-2x)^5$, as claimed.
\end{proof}

\begin{corollary}\label{cor:exactVariance}
For $n\geq1$,
\begin{equation}\label{eq:exactVariance}
 \boxed{
 \Var(X_n)=\frac{[x^n]F_2(x)}{2^{n-1}}
 -\left(\frac{n^2}{3}-\frac{4n}{9}
   +\frac{2}{27}-\frac{2(-1)^n}{27\cdot2^n}\right)^2.}
\end{equation}
In particular, the variance is computable by coefficient extraction
from a fixed rational function.
\end{corollary}
\begin{proof}
The coefficient $[x^n]F_2(x)$ is the sum of the squared statistic
over all compositions of $n$. Dividing by $2^{n-1}$ gives $\E[X_n^2]$.
Subtract the square of the expectation in Corollary~\ref{cor:exactMean}.
\end{proof}

The rationality of $F_1$ and $F_2$ is not restricted to the first
two moments.

\begin{theorem}\label{thm:allMoments}
For every fixed nonnegative integer $r$, the generating function
\[
 F_r(x)=\sum_{n\geq1}\left(\sum_{\alpha\models n}
                                  g_C(\alpha)^r\right)x^n
\]
is rational in $x$.
\end{theorem}
\begin{proof}
The identity $2g_C=J+D$ reduces the problem to products of powers
of $J=n(n-k)$ and of $D=\sum_{i<j}|\alpha_i-\alpha_j|$.
A fixed power $D^b$ expands into a sum over $b$ ordered pairs of
part positions. Every such term involves at most $2b$ distinct
marked positions. For each fixed pattern of coincidences among
these positions, a composition consists of those finitely many
marked positive parts separated by arbitrary compositions. The
latter contribute powers of $C(x,u)$, a rational function.

The remaining contribution is a weighted sum over a fixed number
$h\leq2b$ of positive sizes $a_1,\ldots,a_h$, whose weight is a
product of expressions $|a_i-a_j|$. Split the summation into the
finitely many possible weak orderings of the $a_i$ (including
possible equalities). On each such region the absolute differences
are linear expressions, and the weight is polynomial in the sizes.
Choose the distinct ordered sizes and express their successive
differences as nonnegative integer variables, using positive
differences between distinct sizes. The exponent $a_1+\cdots+a_h$
is then a linear form in these variables. Summing a polynomial
weight times a monomial over a finite product of geometric series
produces a rational function. For example, the needed one-variable
building blocks have the form
\[
 \sum_{d\geq0}d^j x^{ad}
 =\left(x^a\frac{d}{d(x^a)}\right)^j\frac1{1-x^a}
 \qquad(a\geq1).
\]
Finally, powers of $n$ and $k$ are inserted by the operators
$\Theta_x$ and $\Theta_u$ applied a fixed number of times;
differentiation preserves rationality. Every term in the expansion
of $g_C^r$ is therefore rational, and there are only finitely many
such terms for fixed $r$.
\end{proof}

\begin{proposition}\label{prop:nonrational}
The full distribution series $\mathcal H(q,x)$ of
\eqref{eq:fullDistribution} is not rational over $\mathbb Q(q,x)$.
\end{proposition}
\begin{proof}
For the one-part composition $(n)$,
\[
 g_C((n))=\binom n2.
\]
Consequently, for each fixed real $q>1$,
\[
 [x^n]\mathcal H(q,x)\geq q^{\binom n2}.
\]
The $n$th root of the right-hand side tends to infinity, so
$\mathcal H(q,x)$ has radius of convergence zero as a series in $x$.
If $\mathcal H(q,x)$ were rational over $\mathbb Q(q,x)$, its formal
power-series expansion at $x=0$ could be represented by a fraction
whose denominator has nonzero constant coefficient in $\mathbb Q(q)$.
For all but finitely many real specializations $q>1$, that fraction
would specialize to a rational function analytic in a neighborhood
of $x=0$, hence with positive radius of convergence. This
contradicts the coefficient lower bound.
\end{proof}

\section{Asymptotic expectation and variance}\label{sec:asymptotics}

Rational generating functions make the leading asymptotics
particularly transparent. Their coefficient growth is determined
by poles of smallest modulus. In both \eqref{eq:F1} and
\eqref{eq:F2}, the unique dominant pole is $x=1/2$; the poles
at roots of unity have strictly larger modulus. We use the
standard coefficient formula
\begin{equation}\label{eq:poleCoefficient}
 [x^n](1-2x)^{-j}=2^n\binom{n+j-1}{j-1}.
\end{equation}

\begin{theorem}\label{thm:moments}
For a uniform random composition $C_n$, let $X_n=g_C(C_n)$.
As $n\to\infty$,
\begin{align}
 \E[X_n]&=\frac{n^2}{3}-\frac{4n}{9}+O(1),
                         \label{eq:meanAsymptotic}\\
 \Var(X_n)&=\frac{10}{189}n^3+O(n^2).
                         \label{eq:varianceAsymptotic}
\end{align}
For the normalized statistic $\widetilde X_n=\widetilde g_C(C_n)$,
defined when $n\geq2$,
\begin{align}
 \E[\widetilde X_n]
   &=\frac23-\frac{2}{9n}+O(n^{-2}),
                      \label{eq:normalizedMean}\\
 \Var(\widetilde X_n)
   &=\frac{40}{189n}+O(n^{-2}).
                      \label{eq:normalizedVariance}
\end{align}
In particular, $\widetilde X_n\to2/3$ in probability.
\end{theorem}

\begin{proof}
The expectation follows immediately from
\eqref{eq:exactMean}; its final term is exponentially small.
For the second moment, expand $F_2$ about its dominant pole:
\begin{equation}\label{eq:principalPart}
 F_2(x)=\frac{4}{3(1-2x)^5}
        -\frac{256}{63(1-2x)^4}
        +O((1-2x)^{-3})
 \qquad(x\to\tfrac12).
\end{equation}
Here the $O$-term describes the remaining terms of the Laurent
expansion; the exact rational function has no other poles of
modulus $1/2$. Applying \eqref{eq:poleCoefficient} gives
\begin{align*}
 \E[X_n^2]
 &=\frac{[x^n]F_2(x)}{2^{n-1}}\\
 &=\frac83\binom{n+4}{4}
   -\frac{512}{63}\binom{n+3}{3}+O(n^2)\\
 &=\frac{n^4}{9}-\frac{46n^3}{189}+O(n^2).
\end{align*}
On the other hand, \eqref{eq:exactMean} gives
\[
 \E[X_n]^2=\frac{n^4}{9}-\frac{8n^3}{27}+O(n^2).
\]
The $n^4$ terms cancel in the difference. Therefore,
\[
 \Var(X_n)=\left(-\frac{46}{189}+\frac8{27}\right)n^3
           +O(n^2)=\frac{10}{189}n^3+O(n^2).
\]
Finally, $\widetilde X_n=2X_n/[n(n-1)]$ gives
\[
 \E[\widetilde X_n]=\frac{2\E[X_n]}{n(n-1)},\qquad
 \Var(\widetilde X_n)=\frac{4\Var(X_n)}{n^2(n-1)^2}.
\]
Expanding these expressions yields \eqref{eq:normalizedMean} and
\eqref{eq:normalizedVariance}. The mean converges to $2/3$ and
the variance tends to zero, so Chebyshev's inequality proves
convergence in probability.
\end{proof}

\begin{remark}\label{rem:comparison}
For the Gini statistic on uniformly random partitions, the same
normalization tends to $1$ in probability, whereas the composition
model has an interior limit $2/3$, approached at rate $n^{-1}$.
For completeness, the partition limit can be seen directly from the
classical Hardy--Ramanujan asymptotic for the partition numbers
\cite{AndrewsEriksson}. If $\Lambda_n$ is a uniform random partition,
conjugation and the bound
\[
 M_0(\lambda)=\sum_j\binom{\lambda_j^*+1}{2}
 \leq\frac{n(\ell(\lambda)+1)}2
\]
show that it suffices to prove $\ell(\Lambda_n)/n\to0$ in
probability. By conjugation, the number of partitions with at least
$\varepsilon n$ parts equals the number with largest part at least
$\varepsilon n$. Removing that largest part gives, with
$a=\lceil\varepsilon n\rceil$,
\[
 \Prob(\ell(\Lambda_n)\geq\varepsilon n)
 \leq\frac{\sum_{k=a}^n p(n-k)}{p(n)}
 \leq\frac{n p(n-a)}{p(n)}\longrightarrow0
 \qquad(0<\varepsilon<1).
\]
The last limit follows since $\log p(n)$ has leading term
$\pi\sqrt{2n/3}$. Hence the normalized partition Gini statistic
converges to $1$ in probability. The composition model not only
changes this limit, but also has normalized variance of order
$n^{-1}$, as shown above.
\end{remark}

\section{A central limit theorem}\label{sec:CLT}

The rational generating function determines the asymptotic
variance, but does not by itself explain the probabilistic
mechanism leading to a normal law. We give a direct proof using the
natural random-separator construction of a uniform composition.
This also permits a joint limit theorem involving the number of
parts, from which a limit theorem for the ordinary classical Gini
coefficient follows.

\subsection{Random separators and a linearized statistic}

For every $n\geq1$, place separators independently into the $n-1$
gaps between consecutive units, each with probability $1/2$.
Each subset of gaps is equally likely, so the resulting composition
$C_n$ is uniform on $\Comp_n$. Its number of parts satisfies
\begin{equation}\label{eq:numberParts}
 K_n:=K(C_n)\ \overset{d}=\ 1+\operatorname{Bin}(n-1,1/2).
\end{equation}

Let $c_j=c_j(C_n)$ count the parts of size at least $j$, and put
$q_j=2^{-j}$. Define
\begin{equation}\label{eq:Astat}
 A_n=\sum_{j\geq1}q_jc_j,
 \qquad \delta_{j,n}=c_j-nq_j.
\end{equation}
Because $\sum_{j\geq1}q_j^2=1/3$, the column formula
\eqref{eq:columnGini} gives the exact identity
\begin{align}
 g_C(C_n)
 &=\frac12\left(n^2-\sum_j(nq_j+\delta_{j,n})^2\right)
 \notag\\
 &=\frac{n^2}{3}
   -n\left(A_n-\frac n3\right)
   -\frac12\sum_{j\geq1}\delta_{j,n}^2.
 \label{eq:linearization}
\end{align}
The first two terms are a constant and a linear function of the
part counts. We next show that the final, quadratic term is
asymptotically negligible on the scale $n^{3/2}$.

\begin{lemma}\label{lem:quadraticError}
With the preceding notation,
\begin{equation}\label{eq:quadraticBound}
 \E\left[\sum_{j\geq1}\delta_{j,n}^2\right]=O(n).
\end{equation}
Consequently,
$\sum_{j\geq1}\delta_{j,n}^2=o_{\Prob}(n^{3/2})$.
\end{lemma}
\begin{proof}
For $1\leq j\leq n$, the event that a part of size at least
$j$ begins at the first unit has probability $2^{-(j-1)}$.
For a possible beginning at a later unit, the preceding gap must
contain a separator and the following $j-1$ gaps must contain
none, so the probability is $2^{-j}$. Thus
\begin{equation}\label{eq:cMean}
 \E[c_j]=2^{-(j-1)}+(n-j)2^{-j}=(n-j+2)2^{-j},
 \qquad 1\leq j\leq n.
\end{equation}
In particular, $\E[c_j]-nq_j=(2-j)2^{-j}$.

To bound the variance, write $c_j$ as a sum of indicators indexed by
possible beginning positions. Each indicator depends on at most $j$
consecutive separator variables, and hence is independent of all
but at most $2j$ other indicators. Every indicator has expectation
at most $2^{-(j-1)}$. Bounding the covariance of two indicators
in absolute value by a constant times the smaller of their
expectations gives, with a universal constant $C$,
\[
 \Var(c_j)\leq Cnj2^{-j}.
\]
Therefore,
\[
 \E[\delta_{j,n}^2]
 =\Var(c_j)+(\E[c_j]-nq_j)^2
 \leq Cnj2^{-j}+(j-2)^24^{-j}.
\]
The sums $\sum_{j\geq1}j2^{-j}$ and
$\sum_{j\geq1}(j-2)^24^{-j}$ converge. For $j>n$,
$c_j=0$ and $\delta_{j,n}=-n2^{-j}$, whose squared sum is
$O(n^24^{-n})$. Summing proves \eqref{eq:quadraticBound}.
The probabilistic assertion follows from Markov's inequality.
\end{proof}

There is a useful alternative expression for $A_n$. Each part
of size $a$ contributes to $c_j$ exactly when $1\leq j\leq a$.
Hence, with
\begin{equation}\label{eq:hReward}
 h(a)=\sum_{j=1}^a2^{-j}=1-2^{-a},
\end{equation}
we have
\begin{equation}\label{eq:renewalRewardRepresentation}
 A_n=\sum_{i=1}^{K_n}h(\alpha_i),
 \qquad C_n=(\alpha_1,\ldots,\alpha_{K_n}).
\end{equation}
Thus $A_n$ is an \emph{additive} statistic of the parts, even
though $g_C$ itself is quadratic in the counts $c_j$.

\subsection{The renewal-reward central limit theorem}

Couple all the compositions by using an infinite sequence of
independent fair separator choices. The successive complete part
lengths $L_1,L_2,\ldots$ are then independent and geometrically
distributed:
\begin{equation}\label{eq:geometricL}
 \Prob(L_i=a)=2^{-a},\qquad a=1,2,\ldots.
\end{equation}
Let $S_0=0$, $S_r=L_1+\cdots+L_r$, and set
\[
 N_n=\max\{r\geq0:S_r<n\},\qquad R_n=n-S_{N_n}.
\]
The parts of $C_n$ are $L_1,\ldots,L_{N_n},R_n$, and
$K_n=N_n+1$. The last part $R_n$ can be truncated, but is
stochastically bounded: for $1\leq r\leq n$,
\begin{equation}\label{eq:residualTail}
 \Prob(R_n\geq r)=2^{-(r-1)}.
\end{equation}
In particular $R_n=O_{\Prob}(1)$. Also, either by the strong law
of large numbers for $S_r$ or directly from \eqref{eq:numberParts},
\begin{equation}\label{eq:renewalRate}
 \frac{N_n}{n}\longrightarrow\frac12
 \qquad\text{almost surely}.
\end{equation}

The geometric distribution and the reward function in
\eqref{eq:hReward} have elementary moments:
\begin{equation}\label{eq:rewardMoments}
\begin{gathered}
 \E[L]=2,\qquad \E[L^2]=6,\\
 \E[h(L)]=\frac23,\qquad
 \E[h(L)^2]=\frac{10}{21},\qquad
 \E[Lh(L)]=\frac{14}{9}.
\end{gathered}
\end{equation}
For example,
\[
 \E[2^{-L}]=\sum_{a\geq1}4^{-a}=\frac13,
 \qquad
 \E[4^{-L}]=\sum_{a\geq1}8^{-a}=\frac17,
\]
so $\E[h(L)]=1-1/3$ and
$\E[h(L)^2]=1-2/3+1/7$. The other moments follow from
$\sum_{a\geq1}a t^a=t/(1-t)^2$.

Set
\begin{equation}\label{eq:centeredRewards}
 Y_i=h(L_i)-\frac{L_i}{3},\qquad
 Z_i=1-\frac{L_i}{2}.
\end{equation}
Both random variables have mean zero. Direct substitution of
\eqref{eq:rewardMoments} gives
\begin{equation}\label{eq:centeredCov}
 \E[Y_i^2]=\frac{20}{189},\quad
 \E[Z_i^2]=\frac12,\quad
 \E[Y_iZ_i]=\frac29.
\end{equation}
The decompositions at the last complete renewal read
\begin{align}
 A_n-\frac n3
   &=\sum_{i=1}^{N_n}Y_i
       +h(R_n)-\frac{R_n}{3},\label{eq:Acentered}\\
 K_n-\frac n2
   &=\sum_{i=1}^{N_n}Z_i
       +1-\frac{R_n}{2}.\label{eq:Kcentered}
\end{align}
The terminal terms are $O_{\Prob}(1)$ by
\eqref{eq:residualTail}. We shall use the following standard
random-index form of the central limit theorem, including its
short proof for completeness.

\begin{lemma}\label{lem:randomIndexCLT}
Let $V_1,V_2,\ldots$ be independent and identically distributed
centered random vectors in $\mathbb R^d$ with finite covariance
matrix $\Sigma_0$. If nonnegative integer-valued random indices $M_n$
satisfy $M_n/n\to\rho>0$ in probability, then
\[
 \frac1{\sqrt n}\sum_{i=1}^{M_n}V_i
 \xrightarrow{d}\N_d(0,\rho\Sigma_0).
\]
The indices need not be independent of the vectors.
\end{lemma}
\begin{proof}
Fix a vector $v\in\mathbb R^d$ and apply the argument to the
scalar sequence $v^TV_i$. Write $m_n=\lfloor\rho n\rfloor$.
By the ordinary central limit theorem,
$n^{-1/2}\sum_{i=1}^{m_n}v^TV_i$ converges to a centered normal
variable with variance $\rho v^T\Sigma_0v$.

It suffices to show that replacing $m_n$ by $M_n$ changes the sum
by $o_{\Prob}(\sqrt n)$. For any fixed $\varepsilon>0$, the event
$|M_n-m_n|\leq\varepsilon n$ has probability tending to one as
$\varepsilon\downarrow0$ after $n\to\infty$. On this event the
difference is bounded by the maximum of the forward and backward
partial sums of at most $\lceil\varepsilon n\rceil$ centered
variables adjacent to $m_n$. Kolmogorov's maximal inequality
bounds the probability that either such maximum exceeds
$\delta\sqrt n$ by
$O(\varepsilon/\delta^2)$, with a constant depending on
$v^T\Sigma_0v$. First take $n\to\infty$ and then
$\varepsilon\downarrow0$. The difference is $o_{\Prob}(\sqrt n)$.
The scalar assertion follows by Slutsky's theorem; the vector
assertion follows from the Cram\'er--Wold device.
\end{proof}

\begin{theorem}\label{thm:jointCLT}
Let $C_n$ be a uniform random composition of $n$, and put
$X_n=g_C(C_n)$, $\widetilde X_n=X_n/\binom n2$, and
$K_n=K(C_n)$. Then
\begin{equation}\label{eq:jointAK}
 \frac1{\sqrt n}
 \begin{pmatrix}A_n-n/3\\ K_n-n/2\end{pmatrix}
 \xrightarrow{d}
 \N_2\!\left(\begin{pmatrix}0\\0\end{pmatrix},
 \begin{pmatrix}10/189&1/9\\1/9&1/4\end{pmatrix}\right).
\end{equation}
In particular,
\begin{equation}\label{eq:rawCLT}
 \boxed{\frac{X_n-n^2/3}{n^{3/2}}
          \xrightarrow{d}\N\!\left(0,\frac{10}{189}\right),}
\end{equation}
and
\begin{equation}\label{eq:normalizedCLT}
 \boxed{\sqrt n\left(\widetilde X_n-\frac23\right)
          \xrightarrow{d}\N\!\left(0,\frac{40}{189}\right).}
\end{equation}
More precisely, the joint limit of
$\sqrt n(\widetilde X_n-2/3)$ and $(K_n-n/2)/\sqrt n$
is centered bivariate normal with covariance matrix
\begin{equation}\label{eq:jointGK}
 \begin{pmatrix}40/189&-2/9\\-2/9&1/4\end{pmatrix}.
\end{equation}
\end{theorem}
\begin{proof}
Apply Lemma~\ref{lem:randomIndexCLT} to the i.i.d. vectors
$(Y_i,Z_i)$ in \eqref{eq:centeredRewards}, using
$N_n/n\to1/2$ and the covariance values in
\eqref{eq:centeredCov}. The limiting covariance matrix is one
half of the matrix with entries in \eqref{eq:centeredCov}.
The terminal terms in \eqref{eq:Acentered} and
\eqref{eq:Kcentered} vanish after division by $\sqrt n$, so this
proves \eqref{eq:jointAK}.

By Lemma~\ref{lem:quadraticError}, the last term in
\eqref{eq:linearization}, divided by $n^{3/2}$, converges to zero
in probability. Therefore
\[
 \frac{X_n-n^2/3}{n^{3/2}}
 =-\frac{A_n-n/3}{\sqrt n}+o_{\Prob}(1),
\]
which proves \eqref{eq:rawCLT}. Since
$\widetilde X_n=2X_n/[n(n-1)]$, it follows that
\[
 \sqrt n\left(\widetilde X_n-\frac23\right)
 =-2\frac{A_n-n/3}{\sqrt n}+o_{\Prob}(1).
\]
This proves \eqref{eq:normalizedCLT}. Applying the same linear
transformation to the first coordinate of \eqref{eq:jointAK}
produces the covariance matrix \eqref{eq:jointGK}.
\end{proof}

\begin{remark}\label{rem:centerMean}
Centering $X_n$ at its exact expectation, rather than at $n^2/3$,
does not change \eqref{eq:rawCLT}, since
$\E[X_n]-n^2/3=O(n)=o(n^{3/2})$. Similarly,
$\E[\widetilde X_n]-2/3=O(n^{-1})=o(n^{-1/2})$.
The variance constants in the central limit theorem agree with
the independent coefficient asymptotics of
Theorem~\ref{thm:moments}.
\end{remark}

\section{Comparison with the classical Gini coefficient}
\label{sec:classical}

The constant $2/3$ is not the same as the limiting Gini coefficient
of the individual parts of the composition. This distinction is
particularly relevant to the motivating comparison with statistical
measurements of inequality.

\begin{theorem}\label{thm:classicalCLT}
For a uniformly random composition $C_n$, let
$G_n=G_{\mathrm{cl}}(C_n)$, with the convention of
\eqref{eq:classicalGini}. Then
\begin{equation}\label{eq:classicalLimit}
 G_n\xrightarrow{\Prob}\frac13,
 \qquad \E[G_n]\longrightarrow\frac13,
\end{equation}
and
\begin{equation}\label{eq:classicalCLT}
 \boxed{\sqrt n\left(G_n-\frac13\right)
    \xrightarrow{d}\N\!\left(0,\frac{20}{189}\right).}
\end{equation}
\end{theorem}
\begin{proof}
The exact relation \eqref{eq:twoGinis} gives
\[
 G_n=1+\frac{(n-1)\widetilde X_n/n-1}{K_n/n}.
\]
Theorem~\ref{thm:jointCLT} and \eqref{eq:numberParts}
imply $\widetilde X_n\to2/3$ and $K_n/n\to1/2$ in probability.
Therefore
\[
 G_n\xrightarrow{\Prob}1+\frac{2/3-1}{1/2}=\frac13.
\]
The finite-sample coefficient satisfies $0\leq G_n\leq1$.
Convergence in probability and boundedness imply convergence of
expectations to $1/3$.

For the central limit theorem, apply the first-order expansion of
$f(x,y)=1+(x-1)/y$ around $(x,y)=(2/3,1/2)$. Its partial
derivatives at that point are
\[
 f_x(2/3,1/2)=2,\qquad f_y(2/3,1/2)=\frac43.
\]
The factor $(n-1)/n$ multiplying $\widetilde X_n$ changes the
argument by only $O(n^{-1})$, which is negligible on the
$n^{-1/2}$ scale. Thus
\begin{align*}
 \sqrt n\left(G_n-\frac13\right)
 &=2\sqrt n\left(\widetilde X_n-\frac23\right)
   +\frac43\frac{K_n-n/2}{\sqrt n}+o_{\Prob}(1).
\end{align*}
By \eqref{eq:jointGK}, the variance of the limiting normal
variable is
\begin{align*}
 4\left(\frac{40}{189}\right)
 +\frac{16}{9}\left(\frac14\right)
 +2\cdot2\cdot\frac43\left(-\frac29\right)
 &=\frac{160}{189}+\frac49-\frac{32}{27}\\
 &=\frac{20}{189}.
\end{align*}
This proves \eqref{eq:classicalCLT}.
\end{proof}

For additional intuition, the geometric variable $L$ in
\eqref{eq:geometricL} has mean $2$. If $L'$ is an independent copy,
then
\begin{equation}\label{eq:geomDifference}
 \E|L-L'|=\frac43.
\end{equation}
Indeed, the tail formula for nonnegative integer variables gives
\[
 \E[\min(L,L')]
 =\sum_{j\geq1}\Prob(L\geq j)^2
 =\sum_{j\geq1}4^{-(j-1)}=\frac43,
\]
whence
$\E|L-L'|=\E L+\E L'-2\E\min(L,L')=4-8/3=4/3$.
The population Gini coefficient of this geometric distribution is
therefore
\[
 \frac{\E|L-L'|}{2\E L}=\frac{4/3}{4}=\frac13.
\]
This independently identifies the constant in
\eqref{eq:classicalLimit}. Equation~\eqref{eq:twoGinis} then
explains the different limit of the partition-style normalization:
with $K_n/n\to1/2$, the two contributions converge to $1/2$
and $(1/2)(1/3)$, whose sum is $2/3$.

\section{Concluding remarks}\label{sec:conclusion}

Rearranging a composition into a partition gives a natural way to
extend the discrete Gini statistic in \cite{Kopitzke2020}.
The rearrangement-invariant identity
\[
 2g_C(\alpha)=n(n-K(\alpha))
                   +\sum_{i<j}|\alpha_i-\alpha_j|
\]
permits both combinatorial enumeration by distinguished positions
and a probabilistic description based on geometric part lengths.
The two methods are complementary. The rational first- and
second-moment generating functions yield exact formulas and
sharp moment asymptotics; the renewal-reward representation proves
a Gaussian limit theorem and explains the variance constant.

The uniform-partition and uniform-composition models thus have
qualitatively different normalized limits. In the latter,
$\widetilde g_C\to2/3$ with asymptotically normal fluctuations of
order $n^{-1/2}$. The ordinary Gini coefficient of the parts has
a different limit, $1/3$, and also satisfies a central limit theorem.
The two limiting constants are connected by the asymptotic part
density $K_n/n\to1/2$.

Several extensions suggest themselves. One may seek closed forms
for higher moments, investigate large deviations away from the
Gaussian scale, or replace fair separators by biased separators,
which leads to geometric part lengths with a general parameter.
One could also compare the full finite-$n$ distributions under
the uniform composition model with the induced nonuniform measure
on partitions. The rationality theorem suggests that fixed moments
remain algebraically tractable even though their full bivariate
distribution series is not rational.

\bigskip
\noindent\textit{2020 Mathematics Subject Classification.}
Primary 05A05, 60C05; Secondary 05A15, 60F05, 91B82.

\end{document}